\documentclass[11pt,english]{amsart}

\makeatletter
\def\specialsection{\@startsection{section}{1}%
  \z@{\linespacing\@plus\linespacing}{.5\linespacing}%
  {\normalfont}}
\def\section{\@startsection{section}{1}%
  \z@{.7\linespacing\@plus\linespacing}{.5\linespacing}%
  {\normalfont\bfseries\scshape}}
\makeatother

\usepackage{amssymb, amsmath, amsfonts, url, graphicx,verbatim,scrextend}
\usepackage{amscd}

\theoremstyle{plain}
\newtheorem*{thm*}{Theorem}
\theoremstyle{plain}
\newtheorem{thm}{Theorem}[section]
\theoremstyle{definition}

\theoremstyle{plain}

\theoremstyle{plain}

\theoremstyle{plain}
\newtheorem{cor}[thm]{Corollary}
\theoremstyle{remark}

\theoremstyle{remark}
\newtheorem*{acknowledgement*}{Acknowledgement}

\title{Expressing the curvature tensor and connection of a given metric in terms of those of another metric}
\author{Dan Gregorian Fodor}
\date{}

\begin{document}

\maketitle

\pagenumbering{roman}

\pagestyle{myheadings}\markboth{}{}

\pagenumbering{arabic}
\pagestyle{myheadings}

\section*{Introduction}
Given a Riemannian metric $g$ on a manifold, a second metric $m$ can be thought of as a field of symmetric, positive-definite, rank-2 covariant tensors with respect to $g$. We give expressions of $m$'s associated connection, and Riemann curvature tensor  $R_m$, in terms of $R_g$ and certain combinations of covariant derivatives of $m$ (with respect to the Levi-Civita connection associated with $g$). The formulas turn out to be generalizations of the coordinate expressions. Coordinate expression formulas can be recovered from ours by setting $g$ as the Euclidean metric induced by a given coordinate chart. As the covariant derivative induced by $g$ becomes the ordinary partial derivative and the $R_g$ tensor vanishes, the formulas coincide with the well-known coordinate expressions for $m$'s connection and curvature tensor.

\section{Preliminaries}
The coordinate-invariant tensor fields that can be built as polynomials from a given metric $g$ and its partial derivatives have been well studied.These are known as the \textit{metric invariants}.They coincide with the tensor fields built out of $g$, the Riemann tensor $R_{abcd}$ induced by $g$, and its iterated covariant derivatives (which are known as the \textit{curvature invariants})$o_j$ \cite{B1}. Therefore the covariant derivative and Riemann tensor can be used as 'building blocks' for all other invariants. We prove relationships between the curvature and connection of two different metrics, $m$ and $g$, such that all expressions written using covariant derivatives and curvature induced by $m$ admit a rewriting in terms of those induced by $g$.

Since a particular metric tensor-field $m$ uniquely determines both its connection and curvature tensor in a coordinate-invariant way, it is to be expected that these admit an expression from $m$ using the connection and curvature induced by a second metric $g$.

\section{Notational conventions}
As we shall be using covariant derivatives from multiple metrics throughout the paper, we need a way to distinguish them. All formulas will have a first parameter specifying the metric who's connection gives the covariant derivative used in them. They may also have additional parameters signifying the variables used in the formula. All raising and lowering of indices will be made explicit, to avoid confusion with regards to what metric is used. For example:

$F_{e}^{d}(g_{ab}, v^{c}) = v^d_{;e}$ \\This is a formula specifying the covariant derivative of a vector $v^d$ (the second parameter) in terms of the connection induced by the standard first parameter. Where we don't need to specify all the variables, we may shorten it to $v^d_{;e}(g_{ab})$ or $v^d_{;e}(g)$.\\\\
We define the Euclidean metric induced by a coordinate chart to be the metric which, in that specific coordinate chart, has the form: $g_{ab}(p^1,p^2,... p^n) =  \delta_{ab}$ .It is defined only on the domain of the coordinate chart. We shall use $\delta_{ab}$ to denote such a metric. We have  $v^d_{;e}(\delta) = v^d_{,e}$\hspace*{0.0ex}, where $\delta$ is induced by our local chart.
\section{Recovering the connection of a metric}
Given two metrics $m$ and $g$, we are interested in writing the connection of  $m$ in terms of $g$. For an arbitrary vector field $v^{a}$, we know that $v^{a}_{;b}(m)-v^{a}_{;b}(g)$ is tensorial (coordinate independent). From coordinate expansion, we also know that it takes the form $K^{a}_{cb} v^{c}$, where  $K^{a}_{cb}$ is the difference between between the Christoffel symbols of $m$ and $g$ under the local coordinate chart.

We must find an expression for $K^{a}_{cb}$ in terms of the connection for $g$.\\\\
Define $\Gamma^{a}_{bc}(g,m) = \frac{1}{2} m^{an}(m_{nb;c}+m_{nc;b}-m_{bc;n})$. We can recover the Christoffel symbols induced in a local coordinate chart by a metric $g$ as $\Gamma^{a}_{bc}(\delta,g)$. \\
We have $v^{a}_{;b}(m)-v^{a}_{;b}(g)= K^{a}_{cb} v^{c} = (\Gamma^{a}_{cb}(\delta,m) - \Gamma^{a}_{cb}(\delta,g))v^{c}$\\
$\Rightarrow K^{a}_{cb} = \Gamma^{a}_{cb}(\delta,m) - \Gamma^{a}_{cb}(\delta,g)$\\\\
 Choose a coordinate chart $\delta$ that is exponential in a given point $p$, for $g$. At that point, $\Gamma^{a}_{cb}(\delta,g)$ vanishes, leaving $K^{a}_{cb}= \Gamma^{a}_{cb}(\delta,m)$. We know $\Gamma^{a}_{cb}(\delta,m)$ is an expression involving $m$ and the first order partial derivatives of $m$. But, at $p$, first order partial derivatives coincide with the covariant derivatives induced by $g$. For our given $\delta$, $K^{a}_{cb}$ at $p$ takes the form $ \Gamma^{a}_{cb}(\delta,m) = \Gamma^{a}_{cb}(g,m)$. But  $\Gamma^{a}_{cb}(g,m)$  no longer depends on our choice of $\delta$. \\\\
We have given a coordinate-independent expression of the difference between two connections in terms of the connection induced by $g$.
This allows us to write\\ $v^{a}_{;b}(m)-v^{a}_{;b}(g)= \Gamma^{a}_{cb}(g,m) v^{c}$, or \\$v^{a}_{;b}(m) = v^{a}_{;b}(g)+\Gamma^{a}_{cb}(g,m) v^{c}$.\\and gives us the following theorem:

\begin{thm}
Given two metrics $m$ and $g$, and a vector field $v^{a}$, the following relationship holds: $v^{a}_{;b}(m) = v^{a}_{;b}(g)+\Gamma^{a}_{cb}(g,m) v^{c}$.
\end{thm}

Note how writing covariant derivatives induced by $m$ in terms of those induced by $g$ is almost identical to the coordinate expression of the covariant derivative (with the exception of a 'modified' Christoffel symbol. We recover coordinate expressions by setting $g=\delta$ :\\
 $v^{a}_{;b}(m) = v^{a}_{;b}(g)+\Gamma^{a}_{cb}(g,m) v^{c}$ becomes \\ $v^{a}_{;b}(m) = v^{a}_{;b}(\delta)+\Gamma^{a}_{cb}(\delta,m) v^{c}$, which gives \\$v^{a}_{;b} = v^{a}_{,b} +\Gamma^{a}_{cb} v^{c}$ where $\Gamma$ is the Christoffel symbol for $m$. \\
 Knowing that $\Gamma^{a}_{cb}(g,m) = \Gamma^{a}_{cb}(\delta,m) - \Gamma^{a}_{cb}(\delta,g)$, we have:

\begin{cor} \label{i}
Given three metrics $m$, $g$ and $h$, the following relationship holds: $\Gamma^{a}_{cb}(m,g) + \Gamma^{a}_{cb}(g,h) + \Gamma^{a}_{cb}(h,m) =0$.
\end{cor}

In general, we use Christoffel symbols to recover the connection of a given metric, from a coordinate chart. But the above relationship allows us to do the reverse: recover the christoffel symbol of a coordinate chart by applying the covariant connection to its induced Euclidean metric:
$\Gamma^{a}_{cb}(m,g) + \Gamma^{a}_{cb}(g,h) + \Gamma^{a}_{cb}(h,m) =0$ gives \\
$\Gamma^{a}_{cb}(m,m) = 0$, and \\
$\Gamma^{a}_{cb}(g,m) + \Gamma^{a}_{cb}(m,g)  = 0$ \\
In particular: $\Gamma^{a}_{cb}(\delta,g) = -\Gamma^{a}_{cb}(g,\delta)$

The Christoffel symbol of a given chart can be recovered by a formula applying covariant derivatives on its induced Euclidean metric. Therefore, coordinate charts can be studied as geometric objects by means of their induced Euclidean connections.

\section{Recovering the curvature tensor}
Define $R^{l}_{ijk}(g,m)=\Gamma^{l}_{ik;j} - \Gamma^{l}_{ij;k}+ \Gamma^{l}_{js}\Gamma^{s}_{ik} -  \Gamma^{l}_{ks}\Gamma^{s}_{ij}$ \ ,\\ where the free parameters $(g,m)$ are considered passed down to every instance of $\Gamma$ (eg: $\Gamma^{l}_{ik;j}$ is shorthand for $\Gamma^{l}_{ik;j}(g,m)$). The curvature tensor induced by a metric $g$ can be recovered as $R^{l}_{ijk}(\delta ,g)$. \\\\
Since we know how to obtain the connection of $m$ from that of $g$, a viable approach to obtain the curvature tensor is to write it in terms of the connection: \\
$V_{i;jk}(m) - V_{i;kj}(m) =V_l R^l_{ijk}(\delta ,m)$ (Ricci identity).\\\\As we are writing the $m$ connection in terms of $g$, for brevity, any instance of $\Gamma$ will be considered to have default parameters $(g,m)$:\\\\
$V_{i;jk}(m) = (V_{i;j}(m))_{;k}(g) - V_{s;j}(m)\Gamma^{s}_{ik} - V_{i;s}(m)\Gamma^{s}_{jk}=$\\
\hspace*{8ex} $=(V_{i;j}(g)-V_{l}\Gamma^{l}_{ij})_{;k}(g)- (V_{s;j}(g)-V_{l}\Gamma^{l}_{sj})\Gamma^{s}_{ik} - V_{l;s}(m)\Gamma^{s}_{jk}=$\\
\hspace*{8ex} $= V_{i;jk}(g) - V_{l}\Gamma^{l}_{ij;k} -  V_{l;k}\Gamma^{l}_{ij} - V_{s;j}\Gamma^{s}_{ik} + V_{l}\Gamma^{l}_{sj}\Gamma^{s}_{ik} - V_{l;s}(m)\Gamma^{s}_{jk}=$\\
\hspace*{8ex} $=V_{i;jk}(g) +(V_{l}\Gamma^{l}_{sj}\Gamma^{s}_{ik}-V_{l}\Gamma^{l}_{ij;k}) - (V_{l;k}\Gamma^{l}_{ij}+ V_{l;j}\Gamma^{s}_{ik}+V_{l;s}(m)\Gamma^{l}_{jk})$ \\
Antisymmetrising $j$ and $k$, the third group of terms vanishes, and we obtain:\\
$V_{i;jk}(m) - V_{i;kj}(m) = V_{i;jk}(g) - V_{i;kj}(g) + V_l(\Gamma^{l}_{ik;j} - \Gamma^{l}_{ij;k}+ \Gamma^{l}_{js}\Gamma^{s}_{ik} -  \Gamma^{l}_{ks}\Gamma^{s}_{ij})$\\
These group as: \\
$V_l R^{l}_{ijk}(\delta,m)= V_l R^{l}_{ijk}(\delta,g) + V_l(\Gamma^{l}_{ik;j} - \Gamma^{l}_{ij;k}+ \Gamma^{l}_{js}\Gamma^{s}_{ik} -  \Gamma^{l}_{ks}\Gamma^{s}_{ij}) $, giving \\
 $R^{l}_{ijk}(\delta,m)= R^{l}_{ijk}(\delta,g) + (\Gamma^{l}_{ik;j} - \Gamma^{l}_{ij;k}+ \Gamma^{l}_{js}\Gamma^{s}_{ik} -  \Gamma^{l}_{ks}\Gamma^{s}_{ij}) $, or \\
 $R^{l}_{ijk}(\delta,m)= R^{l}_{ijk}(\delta,g) + R^{l}_{ijk}(g,m)$ \\
 Thus we have obtained the curvature tensor of $m$, in terms of the curvature of $g$ and and an application of the connection of $g$ to $m$. \\
 This parallels the process of obtaining a coordinate expression for $R^{l}_{ijk}$ by expanding the anticommutator of the derivatives, with the exception of a remainder as the curvature of $g$, given that the derivatives of the 'base'-metric $g$ do not necessarily commute. The standard formula can be recovered by setting $g = \delta$ .The remainder vanishes, as Euclidean metrics are characterised by null-curvature.

\begin{thm}
Given two metrics $m$ and $g$, the following relationship holds:
$R^{l}_{ijk}(\delta,m)= R^{l}_{ijk}(\delta,g) + R^{l}_{ijk}(g,m)$
\end{thm}

\begin{cor}
Given three metrics $m$ and $g$ and $h$, the following relationship holds:
$R^{l}_{ijk}(m,g) + R^{l}_{ijk}(g,h) + R^{l}_{ijk}(h,m) = 0$
\end{cor}

Analogous relationships are valid for the Ricci tensor, as it is obtained  by contracting the Riemann tensor with $\delta^{j}_{l}$ .An interesting note is that tensors that appear in these kinds of triangular 'sum'-relations are required to be invariant to the resealing of the metric by a constant (this is a necessary but not sufficient condition) . For example, Christoffel symbols of the second kind satisfy a sum-relation (theorem \ref{i}), but Christoffel symbols of the first kind cannot. \\\\

We can also characterize when a metric has flat curvature without appealing to local coordinate charts:
\begin{cor}
Given two metrics $m$ and $g$, $m$ is flat iff \\$R^{l}_{ijk}(\delta,g) + R^{l}_{ijk}(g,m) = 0$
\end{cor}

All the above characterisations make use of $m$ as an object, but use only the connection and curvature of $g$ to derive its geometric properties.

\section{Conclusion}
The above formulas are generalisations of coordinate-expressions, in terms arbitrary (non-Euclidean) background metrics. They may have many applications, for example studying a deformation that flattens a metric, or defining coordinate charts to non-Euclidean domains (using a given Riemann manifold to chart another manifold).

\end{document}